\documentclass[14pt,reqno,a4pape]{amsart}
\usepackage{times,mathrsfs,color,comment}
\usepackage{amssymb,amsmath,bbm}
\usepackage[notref,notcite]{showkeys}

\newtheorem{theorem}{THEOREM}[section]

\newtheorem{lemma}[theorem]{Lemma}

\newtheorem{proposition}[theorem]{Proposition}

\def\wbar{\overline{w}}

\def\nbar{\overline{n}}

\def\jbar{\overline{j}}

\def\hbar{\overline{h}}

\def\Omegabar{\overline{\Omega}}
\def \S{\cal S}
\def\2bar{\overline{2}}

\def\CC{{\rm\kern.24em\vrule
width.02em height1.4ex
depth-.05ex\kern-.26em C}}
\def\QQ{{\rm\kern.24em\vrule width.02em
height1.4ex depth-.05ex\kern-.26em Q}}
\def\RR{{\rm I\kern-.2em R}}
\def\HH{{\rm I\kern-.2em H}}
\def\ZZ{{\rm\kern.26em\vrule width.02em
height0.5ex depth0ex\kern.04em\vrule width.02em
height1.47ex depth-1ex\kern-.34em Z}}
\def\Ibb#1{{\rm I\kern-.23em#1}}
\def\Ib#1{{\rm I\kern-.25em#1}}
\def\k#1{\kern#1em}
\def\vb#1{\vrule width.02em height1.4ex depth-.05ex}
\def\NN{\Ibb N}
\def\11{{\rm\k{.45}\vb0\k{-.142}1}}

\def\ibar{\overline{i}}
\def\epf{\hskip.2in\vrule width.4pt height6.65pt
depth.15pt\vrule
width2.5pt height6.65pt depth-6.25pt\hskip-2.5pt\vrule
width2.5pt
height.25pt depth.15pt\vrule width.4pt
height6.65pt depth.15pt\ }

\def\proof{\noindent {\bf Proof. } }

\def\1bar{\overline{1} }

\def \Omegabar{\overline \Omega}

\def \d{\partial}

\def \S{\widetilde{S}}

\def\zbar{\overline{z}}
\def\dbar{\overline{\partial}}
\def \hbar{\overline{h}}

\def\wbar{\overline{w}}

\def\\GBB{\cal B}

\def\1bar{\overline{1}}

\def\zbar1{\overline{z}^1}

\def\zbar{\overline{z}}

\begin{document}

 \title{\bf   $\dbar$-Estimates  on  the product of   bounded Lipschitz domain }
 
\date{August  27, 2024; 2010 MSC: Primary 32W05, Secondly, 31B10 }

\author{ Song-Ying Li, \    Sujuan Long \ and \ \ Jie Luo}

\maketitle

\centerline{ Dedicated to the memory of Professor Joe Kohn}
\bigskip

\begin{abstract} Let $D$ be a bounded domain in the complex plane with Lipschitz boundary. In the paper, we construct an
integral solution operator $T[f]$ for any $\dbar$ closed $(0,1)$-form $f\in L^p_{(0,1)}(D^n)$ solving
the Cauchy-Riemain equation $\dbar u=f$ on the product domains $D^n$ and obtain the $L^p$-estimates
 for all $1<p\le \infty$.
\end{abstract}

\section{Introduction}

The sup-norm estimates for Cauchy-Riemann equation:
\begin{equation}
\dbar u=f
\end{equation}
on the product domains $\Omega^n$ have received a considerable study recently by many authors.
The research is around the classical problem  posed by Kerzman \cite{Ker} in 1971. After some modification (see \cite{Li}), Kerzman's problem can be stated as follows:
 For any $\dbar$-closed $(0,1)$-form $f\in L^\infty_{(0,1)}(\Omega^n)$,  
 is there  $u\in L^\infty(\Omega^n)$ such that $\dbar u=f$ when $\Omega=D$ is the unit disk in $\CC$? The problem was studied by Henkin \cite{Hen} in 1971,
who proved that if $f\in C^1_{(0,1)}(\overline{D^{n}})$ is $\dbar$-closed,  there is a scalar constant $C$ and $u \in L^\infty(\Omega^n)$ solving $\dbar$-equation (1.1)  such that 
\begin{equation}
\|u\|_{L^\infty(D^n)}\le C\|f\|_{L^\infty_{(0,1)}(D^n)}.
\end{equation}
Notice that $\{f\in C^1_{(0,1)}(\overline{D}^n): \dbar f=0\} $ may not be dense in $\{f\in L^\infty_{(0,1)}(D^n):\dbar f=0\}$. So,  Henkin's
result has not solved the Kerzman's problem. 

Let $A^2(\Omega^n)$ denote the Bergman space consisting of all holomorphic functions $g \in L^2(\Omega^n)$. A solution $u$ of $\dbar$-equation  {(1.1)} is said to be the canonical
solution if $u\perp A^2(\Omega^n)$.  Landucci \cite{Lan} improved  Henkin's result and proved that the estimate (1.2) holds for
 the canonical solution $u$.  Recently, Chen and McNeal \cite{CM} introduced some new
$\tilde{L}^p(D^n)$ and $\tilde{L}^p_{(0,1)}(D^n)$ and obtained $\tilde{L}^p$-estimate for all $1\le p\le \infty$, but their $\tilde{L}^\infty$ is strictly smaller than
$L^\infty(D^n)$. 

In \cite{DPZ}, Dong, Pan and Zhang  proved the canonical solution of $\dbar$-equation (1.1) satisfies the sup-norm estimate (1.2) when $f\in C_{(0,1)}(\Omegabar^n)$
and $\d \Omega\in C^2$. This result greatly improves the previous {work} in \cite{Lan}. However, $\dbar$-closed forms  in $L^\infty_{(0,1)}(\Omega^n)$
may not be approximated by   $\dbar$-closed forms in $C_{(0,1)}(\Omegabar^n)$.
Their result has not solved the Kerzman problem, yet.  Finally, the Kerzman's problem has been settled by Yuan \cite{Yuan} and Li \cite{Li} independently in 2022. They proved that the canonical solution
$u$ of the $\dbar$-equation $(1.1)$ satisfies $\|u\|_{L^p(\Omega^n)}\le C_\Omega^n \|f\|_{L^p_{(0,1)}(\Omega^n)}$ for all $1\le p\le \infty$ with the assumption $\d \Omega\in C^2$ in \cite{Yuan}
and $\d\Omega\in C^{1,\alpha}$ for some $\alpha>0$ in \cite{Li}, respectively. In \cite{Li}, Li  also gives a very beautiful formula for canoncal solution $u$ of $\dbar$-equation (1.1)
which should be very useful for future study (see \cite{Zha}, for example).

    Based on the previous research results, it is very natural to ask the following  question: 

{\it {\bf Question.} Does the sup-norm estimate (1.2) for $\dbar$ hold when $\d \Omega$ is only Lipschitz?}

Let $G(z,w)$ be the Green's function for $\Delta'={\d^2\over \d z\d\zbar}$ on a bounded domain $\Omega\subset \CC$.  in \cite{Li}, Li gives estimates for $G(z,w)$ and its derivatives
when $\d \Omega$ is $C^{1,\alpha}$ for some $\alpha>0$. It is known that the Bergman kernel $K(z,w)$ are related to $G(z,w)$ (see \cite{DPZ}, \cite{Li} and \cite{GL}). In fact
$$
K(z,w)={\d^2 G(z,w)\over \d z\d \wbar},\quad z, w\in \Omega, z\ne w.
$$
From the estimate of $G(z,w)$ and its derivative in \cite{Li}, one has that the Bergman projection $P$ is bounded on $L^p(\Omega)$ for all $1<p<\infty$ and
is bounded from $L^\infty(\Omega)$ to ${\mathcal B}(\Omega)$ (Bloch space) when $\d \Omega\in C^{1,\alpha}$ for some $\alpha>0$.
An example of bounded domain $\Omega_0$ with Lipschitz boundary was constructed by Jerison and Kenig \cite{JK} so that there is a real number $p_1>4$ such that
the Bergman projection $P$ is not bounded on $L^q(\Omega)$ for $p_1'< q< p_1$. Their result indicates that one should consider some solution for $\dbar$ equation (1.1)  instead of the
canonical solution when $\d \Omega$ is Lipschitz. Along this line, we develop some new techniques in this paper so that we are able to answer the above question affirmatively  and prove the following theorem.

\begin{theorem} Let $\Omega$ be a bounded domain in $\CC$ with Lipschitz boundary $\d \Omega$. For any $1<p\le \infty$ and any $\dbar$-closed
$(0,1)$-form $f=\sum_{j=1}^n f_j d\zbar_j\in L^p_{(0,1)}(\Omega^n)$, there is a linear integral operator
\begin{equation}
T[f]=\sum_{j=1}^n T_j[f_j]
\end{equation}
solving $\dbar u=f$ and satisfying the estimate
\begin{equation}
\|T[f]\|_{L^p(\Omega^n)}\le \Big(C_n  C_\Omega \Big)^n \|f\|_{^p_{(0,1)}(\Omega^n)},
\end{equation}
where $C_\Omega$ is constant depending only  on $\|\delta_\Omega(\cdot)\|_\infty+ \| \delta_\Omega (\cdot) \|_{Lip(\Omega)}$, and $\delta_\Omega (z)$ is the distance function from $z$ to $\d \Omega$
and $C_n$ is a positive constant depending only on $n$.
\end{theorem}

For more information about $\dbar$-equations and homotopy formulae, we refer the reader to  the paper of Gong \cite{Gon} and references listed in the current paper.
\medskip

The paper is organized as follows. In Section 2, we give an integral formula for $\dbar$ through mathematical induction. In Section 3, we provide a new method to
transfer the integral formula  solution for  $\dbar$ in Section 2 and get a new integral formula solution for $\dbar$, from which we can get $L^p$ estimate for smooth $f$ with uniform constant given by (1.4) for $1<p\le \infty$\textcolor{red}{.} Finally, in Section 4,
we prove Theorem 1.1.

\section{Integral Formulae Solution I}

%In this section, we suppose
%$$
%f=\sum_{j=1}^n f_j d\zbar_j\in L^p_{(0,1)}(D),
%$$
%\textcolor{red}{is any $\dbar$-closed $(0,1)$-form and $\Omega\subset \CC$  is any domain such that $D\subset \Omega$.} \

Let $D$ and $\Omega$ be two domains in $\CC$ with $D\subset \Omega$.
 By extension theorem stated in \cite{Gon} and \cite{SY} and references therein,
one can extend $f\in W^{k,p}(D)$ so that $f\in W^{k, p}(\Omega)$.
In this section, for any $\dbar$-closed $(0,1)$-form $f\in L^p_{(0,1)}(D)$, we will construct an integral formula $T[f]$ 
solving $\dbar u=f$ on $D$ by introducing a larger domain $\Omega$ which may be useful in future researches. In fact, a homotopy formula for $\dbar$ on $D$ by introducing $\Omega$ 
 with smooth boundary such that $D\subset \Omega$ can help one to solve $\dbar u=f$ on $D$ and
 reduce the smoothness assumption of  $D$.  This technique has been used 
by several authors, such as  Gong \cite{Gon}, Shi and Yao \cite{SY} when $D$ is a strictly pseudoconvex domain
in $\CC^n$ with $C^2$ boundary. They can prove $W^{k,p}(D)$-estimate for $\dbar$ only assuming that $\d D\in C^2$. Using the formula of the canonical solution $T[f]$ for $\dbar$ equation (1,1) constructed in \cite{Li}, in \cite{Zha}, Y. Zhang
proves that if $\d D\in C^\infty$, then the canonical solution  $T[f]\in W^{k,p}(D^n)$ if the $\dbar$-closed $f\in W^{k,p}(D^n)$. Here $k\in \NN$ and $1<p<\infty$. We wish the formula constructed in the section may help to get $W^{k,p}$-estimate $\dbar$
when $\d D$ is only Lipschitz with the help of some extension theorem and technique in \cite{Gon} etc..

Let $g $ be an integrable function on  $\Omega \subset\CC$ and $D\subset \Omega$. Define
\begin{eqnarray}
S_j[g]={1\over 2\pi i} \int_\Omega {1\over z_j-w_j} g (w_j) d\wbar_j\wedge dw_j,\quad 1\le j\le n.
\end{eqnarray}

Let $D$ be a bounded domain in $\Omega$. Define
\begin{eqnarray}
S^D_j[g]={1\over 2\pi i} \int_{\Omega\setminus \overline{D} } {1\over z_j-w_j} g(w_j) d\wbar_j \wedge dw_j,\quad 1\le j\le n.
\end{eqnarray}

\begin{proposition} Let $D$ and $\Omega$ be  bounded domains in $\CC$ with $D\subset \Omega$. Let $f\in C_{(0,1)}^1(\Omega^2)$
be $\dbar$-closed on $D^2$. Define
$$
T^2[f](z)=S_1[f_1]+S_2[f_2]-S_1 S_2\Big[{\d f_2\over \d \zbar_1}\Big]
-S_1^D S_2\Big[ {\d f_1\over \d \zbar_2}-{\d f_2\over \d\zbar_1}\Big].
$$
Then
$
\dbar T^2[f]=f  $ on $ D^2.$

\end{proposition}

\proof  Since ${\d S_j[g]\over \d \zbar_j}=g$ and $S_j^D[g]$ is holomorphic for $z_j\in D$, one has
$$
{\d T^2[f]\over \d \zbar_1}=f_1+{\d \over \d \zbar_1} S_2[ f_2]
-S_2[{\d f_2\over \d \zbar_1}]=f_1,\quad z_1\in D.
$$
For any $z\in D^2$, since $f$ is $\dbar$-closed on $D^2$, one has
$$
{\d f_2\over \d \zbar_1}-{\d f_1\over \d \zbar_2}=0.
$$
Thus, for $z_2\in D$, one has
\begin{eqnarray*}
{\d T^2[f]\over \d \zbar_2}&=& f_2+{\d \over \d \zbar_2} S_1[f_1]
-S_1[{\d  f_2\over \d \zbar_1}] -S_1^D[ {\d f_1 \over \d \zbar_2}-{\d f_2\over \d \zbar_1}] \\
&=&f_2+ 
S_1[  {\d f_1\over \d \zbar_2}-{\d f_2\over \d \zbar_1}]-S_1^D[ {\d f_1 \over \d \zbar_2}-{\d f_2\over \d \zbar_1}] \\
&=& f_2+S_1^D[  {\d f_1\over \d \zbar_2}-{\d f_2\over \d \zbar_1}]-S_1^D[ {\d f_1 \over \d \zbar_2}-{\d f_2\over \d \zbar_1}]\\
&=& f_2.
\end{eqnarray*}
Therefore, the proof of the proposition is complete.\epf

Let $f\in C^1_{(0,1)}(\Omega^n)$  be  $\dbar$-closed on $D^n$.  Write
\begin{eqnarray}
f=f^{n-1}+f_n d\zbar_n,\quad f^{n-1}=\sum_{j=1}^{n-1} f_j d\zbar_j.
\end{eqnarray}
Assume that $T^{n-1}[f^{n-1}]$ has been constructed such that
$$
{\d T^{n-1}[f] \over \d \zbar_j}=f_j, \ \hbox{ on}\quad D^{n-1},\quad 1\le j\le n-1.
$$
Define a $(0,1)$-form in $z_1,\cdots, z_{n-1}\in \Omega$ as follows:
\begin{eqnarray}
{\mathcal D}_{n-1} [f]={\d f^{n-1}\over \d \zbar_n}-\sum_{j=1}^{n-1}{\d f_n \over \d \zbar_j} d\zbar_j={\d f^{n-1}\over \d \zbar_n}-\dbar' f_n,
\end{eqnarray}
%Notice that
%$$
%S^D T^{n-1}[{\mathcal D}_{n-1} f]
%=-S^D[T^{n-1}[\sum_{j=1}^{n-1}{\d f_n\over \d \zbar_j } d\zbar_j]
%+S^D [T^{n-1}[{\d f^{n-1}\over \d\zbar_n}]
%$$
 where
$$
\dbar' f_n=\sum_{j=1}^{n-1}{\d f_n\over \d \zbar_j} d\zbar_j
$$
is $\dbar$-closed in $(z_1,\cdots, z_{n-1})\in \Omega^{n-1}$ for any $z_n\in \Omega$. Then
$$
{\d T^{n-1}[\dbar' f_n ] \over \d \zbar_j}={\d f_n\over \d \zbar_j},\quad 1\le j\le n-1.
$$
\begin{proposition} With the notations above. If $f\in C^1_{(0,1)}(\Omega^n)\cap L^1_{(0,1)}(\Omega)$ is $\dbar$-closed on $D^n$ and if
\begin{equation} \label{T^n}
T^n[f]=S_n[f_n]+T^{n-1}[f^{n-1}]- S_n T^{n-1}[{\d f^{n-1}\over \d \zbar_n}] +S_n^D T^{n-1}[{\mathcal D}_{n-1}(f)],
\end{equation}
then $\dbar T^n [f]=f$ on $D^n$.
\end{proposition}

\proof  By the definition of $T^n[f]$ given by (\ref{T^n}), for any $z_n\in D$, one has
%\begin{eqnarray}
%T^n [f]&=&S_n[f_n]+T^{n-1}[f^{n-1}]- S_n^D T^{n-1}[\dbar' f_n]- (S_n-S_n^D)  T^{n-1}[{\d f^{n-1}\over \d \zbar_n}]\nonumber\\
%&=&S_n[f_n]+T^{n-1}[f^{n-1}]- S_n T^{n-1}[{\d f^{n-1}\over \d \zbar_n}] +S_n^D T^{n-1}[{\mathcal D}_{n-1}(f)] 
%\end{eqnarray}
%Then
$$
{\d T^n [f]\over \d \zbar_n} =f_n+T^{n-1}[{\d f^{n-1}\over \d\zbar_n}]-T^{n-1}[{\d f^{n-1}\over \d \zbar_n}]=f_n.
$$
Notice that
$$
S^D_n [{\mathcal D}_{n-1}]=S^D_n [{\d f^{n-1}\over \d \zbar_n} -\dbar' f_n]=S_n^D[{\d f^{n-1}\over \d \zbar_n}] -S_n^D[\dbar' f_n]
$$
and (\ref{T^n}), one has 
$$
T^n[f]=S_n[f_n]+T^{n-1}[f^{n-1}]- (S_n-S^D)  T^{n-1}[{\d f^{n-1}\over \d \zbar_n}] -S_n^D T^{n-1}[\dbar' f_n].
$$
For $1\le j\le n-1$ and $z\in D^n$, since $\dbar' f_n(w)$ is $\dbar$-closed on $\Omega^{n-1}$ for any $w_n \in \Omega$, one has
\begin{eqnarray*}
{\d T^n [f]\over \d \zbar_j}
&=&S_n\Big[{\d f_n \over \d \zbar_j}\Big]+f_j -(S_n-S_n^D )[{\d f_j \over \d \zbar_n}]-S_n^D [{\d f_n \over \d\zbar_j }]\\
&=&S_n\Big[{\d f_n \over \d \zbar_j}\Big]+f_j -(S_n-S_n^D) [{\d f_n \over \d \zbar_j}]-S_n^D [{\d f_n\over \d\zbar_j}]\\
&=&f_j.
\end{eqnarray*}
The proof of theorem is complete.\epf

\section{Integral Formula Solution II}

In this section, we will transform the formula solution in Section 2 to a new formula solution for $\dbar$ which can be used to get $L^p(D^n)$-estimate
 when $\d D$ is Lipschitz for all $1< p\le \infty.$ We will use {\bf the formula in Section 2 with $D=\Omega$}. In order to do this, we need to develop
some new technique.

For any domain $\Omega$  in $\CC$ with $\Omega \ne \CC$, we define
\begin{equation}
\delta_\Omega (z)=\hbox{dist}(z, \d \Omega)=\inf\{|z-w|: w\in \d \Omega\}
\end{equation}

We say that a bounded domain $\Omega \subset \CC$ is Lipschitz if 
\begin{equation}
\| \nabla \delta_{\Omega}\|_\infty <\infty \iff \delta_\Omega \in \hbox{Lip}(\Omega).
\end{equation}

\medskip
To construct the formula, we borrow some ideas from \cite{Li}. Before we do that, we  introduce some notations here.

For any $i\ne j$, we define

\begin{equation}
\tau_{i,j}(z,w)=|z_i-w_i|^2 \delta(w_i)+|z_j-w_{j}|^{2}\delta(w_j),
\end{equation}

\begin{equation}
A^{i}_{i,j}(z,w)={\d\over\d \wbar_{i}}\Big({(\zbar_{i}-\wbar_{i})\delta(w_{i})\over \tau_{i,j}(z,w)}\Big),\quad A^{j}_{i,j}(z,w)={\d\over \d \wbar_j} \Big( {(\zbar_j-\wbar_j)\delta(w_j)\over \tau_{i,j}(z,w)}\Big)
\end{equation}
and
\begin{equation}
S^i_{i,j}(z,w)={1\over z_i-w_i} A^j_{i,j}(z,w),\quad S^{j}_{i,j}(z,w)={1\over z_j-w_j} A^i_{i,j}(z,w).
\end{equation}
Define %$S^j_j=S_j$, 
\begin{equation}
S^{k}_{i,j} [f_k]= \int_{\Omega^2} S^{k}_{i,j}(z, w) f_k (w){ dv(w)\over \pi^2}, \quad k=i, j.
\end{equation}
%and
%\begin{equation}
%S^{j}_{i,j} [f_j]= \int_{D^2} S^{j}_{i,j}(z, w) f_j (w)  {dv(w) \over \pi^2}.
%\end{equation}
\begin{theorem}\label{3.1} Let $\Omega$ be a bounded domain in $\CC$ with $\delta_\Omega(\cdot)\in \hbox{Lip}(\Omega)$. 
Let $f\in C^1_{(0,1)}(\Omegabar^{n} )$ be $\dbar$-closed. Then there are linear integral operators
\begin{equation}\label{Tj}
T_j[f_j]=S_j[f_j]+\sum_{k=1}^{n-1}\sum_{|I|=k} c_IS^j_J[f_j],\quad S^j_J[f_j]=\int_{\Omega^{k+1}} S^j_J(z, w) f_j(w) {dv(w)\over \pi^{k+1}}
\end{equation}
with $J=\{j, I\} $, $I=\{i_1,\cdots, i_k\}$, $j\not\in I$ and
the linear integral operator
\begin{equation}\label{T}
T^{n}[f]=\sum_{j=1}^n T^{n}_j [f_j]
\end{equation}
satisfying that
$
\dbar T^{n}[f]=f
$ on $\Omega^n$. Moreover,  for each $I=\{i_1,\cdots, i_k\}$, we can write
$$
I=I_1\cup I_2\cdots\cup I_\ell \hbox{ with }|I_\alpha |\ge 1 \  \hbox{ and }  \quad i_\alpha ^*\in I_{\alpha-1} \hbox{ is some element  and }\ i_1^*=j.
$$
Here $I_{\alpha}\cap I_{\beta}=\emptyset$ if $\alpha\neq\beta$ and $\alpha=1,2,\cdots,\ell$.
%Let 
%$$
%\epsilon_{\alpha-1}=\sum_{i \in I_\alpha } {1\over p_i}, \quad p_i>2
%$$
Then
%\begin{equation} \label{S}
%|S^j_J (z,w)|\le  {2^{k-1}\over |w_j-z_j|}  \prod_{\alpha=1}^\ell \prod_{i\in I_\alpha} {(\delta(w_i)+c_0|w_i-z_i|)\over \tau_{i_\ell^* , i}} 
%\end{equation}
\begin{equation} \label{S}
|S^j_J (z,w)|\le  {2^{k-1}\over |w_j-z_j|}  \prod_{s=1}^\ell \prod_{i\in I_s} {(\delta(w_i)+c_0|w_i-z_i|)\over \tau_{i_{s}^* , i}(z,w)}
\end{equation}
and  
\begin{equation}\label{S1}
\Big|{\d S_J^j \over \d \zbar_j}(z,w) \Big|\le   2^{k-1} \max_{i\in I_1} {\delta(w_j)+c_0|w_j-z_j| \over \tau_{j, i}(z,w) } \prod_{\alpha=1}^\ell  \prod_{i\in I_\alpha }{(\delta(w_i)+c_0|w_i-z_i|)\over \tau_{i_{\alpha}^*, i}(z,w)},
\end{equation}
where $c_0=\max\{\|\nabla \delta\|_\infty, 1\}$.
\end{theorem}

\proof We use formulae in Section 2 with $D=\Omega$. When $n=2$, we have
$$
T^2[f]=S_1[f_1]+S_2[f_2]-S_1 S_2[{\d f_2\over \d \zbar_1}].
$$
Since $f$ is $\dbar$-closed, one has
\begin{eqnarray*}
S_1 S_2[{\d f_2\over \d \zbar_1}]
&=&{1\over \pi^2} \int_{\Omega^2}  {1\over (z_1-w_1)(z_2-w_2)} {\d f_2\over \d \wbar_1} d A(w_1) d A(w_2)\\
&=&{1\over \pi^2}\int_{\Omega^2}  {1\over (z_1-w_1)(z_2-w_2)} {|w_1-z_1|^2 \delta(w_1)\over \tau_{1,2}}{\d f_2\over \d \wbar_1} d A(w_1) dA(w_2)\\
&&+{1\over \pi ^2} \int_{\Omega^2}  {1\over (z_1-w_1)(z_2-w_2)} {|w_2-z_2|^2 \delta(w_2)\over \tau_{1,2}}{\d f_1\over \d \wbar_2} d(w_1) d A(w_2)\\
&=&-S^2_{1,2}[f_2]-S^1_{1,2}[f_1].
\end{eqnarray*}
Therefore,
\begin{eqnarray}
T^2[f]=S_1[f_1]+S^1_{1,2}[f_1]+S_2[f_2]+S^2_{1,2}[f_2].
\end{eqnarray}
For any $i\ne j$, since
\begin{eqnarray*}
\lefteqn{A^j_{i,j}(z,w)}\\
&=&{\d \over \d \wbar_j} \Big({(\zbar_j-\wbar_j)\delta(w_j)\over \tau_{i,j} (z,w)}\Big)\\
&=&-{\delta(w_j) \over \tau_{i,j}(z,w)}
+{(\zbar_j-\wbar_j) \over \tau_{i,j}(z,w) } \d_{\jbar} \delta(w_j)+{|w_j-z_j|^2\delta(w_j)\over \tau_{i,j}(z,w)^2}(\delta(w_j) +(\wbar_j-\zbar_j)\d_{\jbar} \delta(w_j))\\
&=&
-{|z_i-w_i|^2 \delta(w_i)   \over \tau_{i,j}(z,w)^2 } \Big(\delta(w_j)+(\wbar_j-\zbar_j) \d_{\jbar} \delta(w_j)\Big),
\end{eqnarray*}
one has
\begin{eqnarray}\label{s1}
S^i_{i,j}(z,w)&=& { A^j_{i,j}\over z_i-w_i}={ (\wbar_i-\zbar_i) \delta(w_i)   \over \tau_{i,j}(z,w) } {\delta(w_j)+(\wbar_j-\zbar_j) \d_{\jbar} \delta(w_j) \over \tau_{i,j}} 
%&=&{(\wbar_i-\zbar_i) \delta(w_i) \delta(w_j)  \over \tau_{i,j}(z,w)^2 } +{(\wbar_j-\zbar_j) \d_{\jbar} \delta(w_j)(\wbar_i-\zbar_i) \delta(w_i)   \over \tau_{i,j}(z,w)^2 }
\end{eqnarray}
and
\begin{eqnarray*}
\lefteqn{{\d \over \d \wbar_i} S^i_{i,j}(z,w)}\\
&=&{(\delta(w_i)+(\wbar_i-\zbar_i)\d_{\ibar} \delta(w_i))(\delta(w_j)+(\wbar_j-\zbar_j)\d_{\jbar} \delta(w_j))
\over \tau_{i,j}^2}\\
&&-2{ |\wbar_i-\zbar_i|^2\delta(w_i) (\delta(w_j)+(\wbar_j-\zbar_j)\d_{\jbar} \delta(w_j))(\delta(w_i)+(\wbar_i-\zbar_i) \d_{\ibar} \delta(w_i))
 \over \tau_{i,j}(z,w)^3 }\\
 &=&{( |w_j-z_j|^2\delta(w_j)-|w_i-z_i|^2\delta(w_i) )(\delta(w_j)+(\wbar_j-\zbar_j)\d_{\jbar} \delta(w_j))(\delta(w_i)+(\wbar_i-\zbar_i) \d_{\ibar} \delta(w_i))
 \over \tau_{i,j}(z,w)^3 }.
\end{eqnarray*}
Therefore, %with $c_0=\max\{\|\nabla \delta\|_{L^\infty(\Omega)}, \hbox{dia}(\Omega), 1\}$,
\begin{eqnarray}
|S^i_{i,j}(z,w) |
&\le &{1\over |w_i-z_i|} {\delta(w_j)+c_0|w_j-z_j|\over \tau_{i,j}(z,w)}
%&\le &{c_0 \delta(w_i)^{-1/\ell} \over |w_i-z_i| ^{1+2/\ell} }  \Big({\delta(w_j)^{1-1/\ell'}\over |w_j-z_j|^{2/\ell'}}+{ |w_j-z_j|^{1-2/\ell' } \over \delta(w_j)^{1/\ell'}}\Big).
\end{eqnarray}
%Let
%\begin{equation}
%\epsilon={1\over p_i}+{1\over p_j}.
%\end{equation}
and
\begin{eqnarray} \label{ds1}
\Big|{\d \over \d \wbar_i} S^i_{i,j}(z,w)\Big|
&\le & {(\delta(w_i)+c_0|w_i-z_i|)\over \tau_{i,j}(z,w) }{(\delta(w_j)+c_0 |w_j-z_j|)\over \tau_{i,j}(z,w)}.
\end{eqnarray}
By the symmetry, this gives the proof of theorem when $n=2$. 
\medskip

When $n>2$ and $f=\sum_{j=1}^n f_j d\zbar _j$ is  $\dbar$-closed, we write
\begin{eqnarray}
f=f^{n-1}+f_n d \zbar_n.
\end{eqnarray}
Assume that
we have constructed
$$
T^{n-1} [f^{n-1}]=\sum_{j=1}^{n-1} T^{n-1}_j [f_j]
$$
with $\dbar' T^{n-1} [f^{n-1}]=f^{n-1}$ and
$$
T^{n-1}_j[f_j]=\sum_{k=1}^{n-1} \sum_{|I|=k} c_I S^j_{j, I}[f_j]
$$
\begin{eqnarray}
S^j_{j, I}[f_j]=\int_{\Omega^k} S^j_{j, I}(z,w) f_j(w) dv_I(w),
\end{eqnarray}
where
$I=\{ i_1,\cdots, i_k \} \subset \{1,2,\cdots, n-1\}$ and $j\not\in I.$
%$$
%S^j_I(z,w)={1\over z_j-w_j} \tilde{S}^j_I (z,w).
%$$
For $I$ as the above and $j\not\in I$ and $j\le n-1$, 
define
\begin{eqnarray} \label{ss}
S^j_{J, n}(z,w)=A_{j, n}^n (z,w) S^j_J(z,w),\quad \S^{j}_{J}(z,w)={S^{j}_{J}(z,w)\over z_{n}-w_{n}}
\end{eqnarray}
and
\begin{eqnarray}\label{ss1}
S^n_{J, n}(z,w)= {\d \over \d \wbar_j} \Big({|w_j-z_j|^2 \delta(w_j) \over  \tau_{j, n} }\S^j_J(z,w)\Big). %\quad \S^{j}_{J}(z,w)={S^{j}_{J}(z,w)\over z_{n}-w_{n}}.
\end{eqnarray}
%where $\S^{j}_{J}(z,w)={S^{j}_{J}(z,w)\over z_{n}-w_{n}}$.
%\begin{eqnarray}\label{ss1}
%S^n_{I, n}(z,w)= {\d \over \d \wbar_j} \Big({(\zbar_j-\wbar_j) \delta(w_j) \over  \tau_{j, n} }\tilde{S}^j_I(z,w)\Big).
%\end{eqnarray}

By (\ref{T^n}), one has
\begin{eqnarray}\label{Tn}
T^n[f]=S_n[f_n]+T^{n-1}[f^{n-1}]+S_n T^{n-1}\Big[{\d f^{n-1}\over \d \zbar_n}\Big].
\end{eqnarray}
Since $f$ is $\dbar$-closed, one has ${\d f_j(w) \over \d\wbar_n}={\d f_n(w) \over \d \wbar_j}$ and
\begin{eqnarray}\label{SS}
\lefteqn{S_n S^j_J [{\d f_j\over \d \wbar_n}] } \nonumber \\
&=&{1\over \pi^{k+2}} \int_{\Omega^{k+2}} {1\over z_n-w_n} S^j_J(z,w){ \d f_j\over \d\wbar_{n}}(w) dv_J dA(w_n) \nonumber \\
&=&-{1\over \pi^{k+2}}\int_{\Omega^{k+2}} {\d \over \d \wbar_n} \Big({(\zbar_n-\wbar_n) \delta(w_n) \over \tau_{j, n} }S^j_J(z,w)\Big) f_j(w) dv_J dA(w_n) \nonumber \\
&&-{1\over \pi^{k+2}}\int_{\Omega^{k+2}} {\d \over \d \wbar_j} \Big({|z_j-w_j|^2 \delta(w_j) \over (z_n-w_n) \tau_{j, n} }S^j_J(z,w)\Big) f_n(w) dv_J dA(w_n) \nonumber \\
&=&-{1\over \pi^{k+2}}\int_{\Omega^{k+2}} A^n_{j, n}(z,w) S^j_J(z,w) f_j(w) dv_J dA(w_n)\nonumber \\
&&-{1\over \pi^{k+2}}\int_{\Omega^{k+2}}  {\d \over \d \wbar_j} \Big({(\zbar_j-\wbar_j) \delta(w_j) \over  \tau_{j, n} }\tilde{S}^j_J(z,w)\Big) f_n(w) dv_{J, n}
\nonumber \\
&=&-{1\over \pi^{k+2} } \int_{\Omega^{k+2}} S^j_{J, n} (z,w)f_j(w) dv_{J, n}
-{1\over \pi^{k+2}} \int_{\Omega^{k+2}} S^n_{J, n} (z,w)f_n(w) dv_{J, n}\nonumber\\
&=&-S^j_{J, n}[f_j]-S^n_{J, n}[f^n].
\end{eqnarray}
Therefore, by the formula for $T^{n-1}[f]$ and $T^n[f]$ given by (\ref{Tn}),
$$
T^n[f]=\sum_{j=1}^n T^n_j[f_j],\quad T_j^n[f_j]=\sum_{k=1}^{n-1} \sum_{|I|=k} c_I S^j_J [f_j]
$$
and
$
\dbar T^n[f]=f.
$
\medskip

Now we assume that (\ref{S}) and (\ref{S1}) hold for $S_J^j$. Next we will prove that they are true for $S^j_{J, n}$ and $S^n_{J, n}$ and the proof of  the theorem is complete
by the Principle of Mathematical Induction.

By (\ref{ss}), (\ref{ss1}) and (\ref{SS}), one has that
$$
S^j_{J,  n}(z, w)=A^n_{j, n} S^j_J(z,w) =-{|w_j-z_j|^2 \delta(w_j)(\delta(w_n)+(\wbar_n-\zbar_n) \d_{\nbar} \delta(w_n))\over \tau_{j,n}^2(z,w)} S_J^j(z,w).
$$
Thus, since $S^j_J$ satisfies (\ref{S}), one has
\begin{eqnarray*}
|S^j_{J, n}(z,w)|&\le&  {\delta (w_n)+c_0|w_n-z_n|\over \tau_{j,n}(z,w) } {2^{k-1}\over |w_j-z_j|}\prod_{s=1}^\ell \prod_{i\in I_s} {\delta(w_{i})+c_0|w_{i}-z_{i}|\over \tau_{i_s^*, i} }\\
&=& {2^{k-1} \over |w_j-z_j| } \prod_{s=1}^{\ell} \prod_{i\in I_s'} {\delta(w_{i})+c_0|w_{i}-z_{i}|\over \tau_{i_s^*, i} }\\
\end{eqnarray*}
with $I'_1=I_1\cup\{n\}$ and $I_s'=I_s$ if $2\le s\le \ell$. This implies that  (\ref{S}) holds for $|I|=k+1$. 

Since
$$
A^n_{j,n}=(z_j-w_j)S^j_{j,n}, \quad S^j_{J, n}=A^n_{j, n} S^j_J=(z_j-w_j)S^j_{j,n} S^j_{J}
$$
and
$$
{\d S^j_{j, n}\over \d \wbar_j}={\delta(w_j)+c_0|w_j-z_j|\over \tau_{j, n} } {\delta(w_n)+ c_0|w_n-z_n| \over \tau_{j, n} },
$$
one has with $I_1'=I_1\cup\{n\}$ and $I_s'=I_s$ if $2\le s\le \ell$ that
\begin{eqnarray*}
\lefteqn{\Big |{\d S^j_{J, n}\over \d \wbar_j} \Big |}\\
%&=& \Big |{\d  (A^n_{j,n} S^j_{J} )\over \d \wbar_j} \Big |\\
 &=& \Big |{\d S_{j, n}^j\over \d \wbar_j } (w_j-z_j)S^j_{J} (z,w)+S^j_{j, n}
(w_j-z_j) {\d S^j_J \over \d\wbar_j }\Big|\\
&\le& {\delta(w_j)+c_0|w_j-z_j|\over \tau_{j, n} } {\delta(w_n)+ c_0|w_n-z_n| \over \tau_{j, n} } |(w_j-z_j)S^j_J|+|S^j_{j, n}| |w_j-z_j| |{\d S^j_J\over \d \wbar_j}|\\
&\le& {\delta(w_j)+c_0|w_j-z_j|\over \tau_{j, n} } {\delta(w_n)+ c_0|w_n-z_n| \over \tau_{j, n} } 2^{k-1}
\prod_{s=1}^\ell \prod_{i\in I_s}  {\delta(w_i)+c_0|w_i-z_i|\over \tau_{i_s^*, i}}\\
&&+2^{k-1}{|w_j-z_j|^2 \delta(w_j) (\delta(w_n)+c_0|w_n-z_n|)  \over \tau_{j, n}^2}
\max_{i\in I_1}  {\delta(w_j) +c_0 |w_j -z_j| \over \tau_{j, i}} 
\prod_{s=1}^\ell\prod_{i\in T_s}  {\delta(w_i)+c_0|w_i-z_i |\over \tau_{i_s^*, i} }\\
&\le& 2^{k-1} {\delta(w_j)+c_0|w_j-z_j|\over \tau_{j, n} }  
\prod_{s=1}^{k} \prod_{i\in I_s'}{\delta(w_i)+c_0|w_{i}-z_{i}|\over \tau_{i_s^*, i } }\\
&&+2^{k-1}
  \max_{i\in I_1'}{\delta(w_j) +c_0 |w_j-z_j| \over \tau_{j, i}} 
\prod_{s=1}^{\ell} \prod_{i\in I_s'} {\delta(w_i)+c_0 |w_i-z_i|\over \tau_{i_s^*, i }}\\
&\le & 2^{k}\max_{i\in I_1'}
   {\delta(w_j) +c_0 |w_j-z_j| \over \tau_{j, i}} 
\prod_{s=1}^{\ell} \prod_{i\in I_s'} {\delta(w_i)+c_0 |w_i-z_i|\over \tau_{i_s^*, i }}.
\end{eqnarray*}
Therefore, (\ref{S1}) holds for $S^j_{J, n}$. 
Next we study $S^n_{J, n}$. Since
\begin{eqnarray*}
 S^n_{J, n} &=&{1\over z_n-w_n}{\d \over \d \wbar_j}\Big({|w_j-z_j|^2 \delta(w_j) \over \tau_{j, n}} S^j_J(z, w)\Big)\\
&=&{1\over z_n-w_n}
\Big( A^j_{j,n}(z,w) (w_j-z_j)S^j_J(z,w) +{|w_j-z_j|^{2} \delta(w_j) \over  \tau_{j, n} }  {\d  S^j_J(z,w) \over \d \wbar_j}\Big) \\
&=& S_{j, n}^n (w_j-z_j) S^j_{J}+{|w_j-z_j|^2 \delta(w_j)\over (z_n-w_n) \tau_{j, n} }{\d S_J^j\over \d \wbar_j}.
\end{eqnarray*}
Thus, by (\ref{S}) and (\ref{S1}), one has
\begin{eqnarray*}
| S^n_{J, n}| %&=&{1\over z_n-w_n}{\d \over \d \wbar_j}\Big({|w_j-z_j|^2 \delta(w_j) \over \tau_{j, n}} S^j_J(z, w)\Big)\\
%&&=&{1\over z_n-w_n}
%\Big( A^j_{j,n}(z,w) (w_j-z_j)S^j_J(z,w) +{(\zbar_j-\wbar_j) \delta(w_j) \over  \tau_{j, n} } (w_j-z_j) {\d  S^j_J(z,w) \over \d \wbar_j}\Big) \\
&\le & |S_{j, n}^n (w_j-z_j) S^j_{J}|+|{|w_j-z_j|^2 \delta(w_j)\over (z_n-w_n) \tau_{j, n} }{\d S_J^j\over \d \wbar_j}|\\
&\le& {2^{k-1} \over |w_n-z_n| }{\delta(w_j)+c_0|w_j-z_j|\over \tau_{j, n}}  \prod_{s=1}^\ell \prod_{i\in I_s}  {\delta(w_{i})+c_0|w_{i}-z_{i}|\over \tau_{i_s^*, i}(z,w)}\\
&&+{|w_j-z_j|^2 \delta(w_j)\over |z_n-w_n| \tau_{j, n} } 2^{k-1} \max_{i\in I_1} {(\delta(w_j)+c_0|w_j-z_j|)\over \tau_{j, i}} \prod_{s=1}^\ell \prod_{i\in I_s}  {\delta(w_{i})+c_0|w_{i}-z_{i}|\over \tau_{i_s^*, i}(z,w)}\\
&\le& {2^{k} \over |w_n-z_n| }\prod_{s=1}^\ell \prod_{i\in I_s'}  {\delta(w_{i})+c_0|w_{i}-z_{i}|\over \tau_{i_s^*, i}(z,w)}\\
\end{eqnarray*}
with $I_1'=I_1\cup\{j\}$ and $I_s'=I_s$ for $2\le s\le \ell$. This implies  that (\ref{S}) holds with $|I|=k+1$.  

Next, we compute and estimate
\begin{eqnarray*}
\lefteqn{\Big|{\d S^n_{J, n}\over \d \wbar_n}\Big|}\\
&\le& \Big|{\d S_{j, n}^n \over \d \wbar_n}  (w_j-z_j) S^j_J \Big|\\
&&+\Big|{|w_j-z_j|^2 \delta(w_j) (\delta(w_n)+(\wbar_n-\zbar_n)\d_{\nbar} \delta(w_n))\over \tau_{j, n}^2 }{\d S_J^j\over \d \wbar_j}\Big|\\
&\le& {(\delta(w_j)+c_0|w_j-z_j|) (\delta(w_n)+c_0|w_n-z_n|)\over \tau_{j, n}^2}
2^{k-1} \prod_{s=1}^\ell \prod_{i\in I_s}  {\delta(w_{i})+c_0|w_{i}-z_{i}|\over \tau_{i_s^*, i}(z,w)}\\
&&+{|w_j-z_j|^2 \delta(w_j) (\delta(w_n)+c_0|w_n-z_n|)\over \tau_{j, n}^2 } 2^{k-1} \max_{i\in I_1}{\delta(w_j)+c_0|w_j-z_j|\over \tau_{j, i}} \prod_{s=1}^\ell \prod_{i\in I_s}  {\delta(w_{i})+c_0|w_{i}-z_{i}|\over \tau_{i_s^*, i}(z,w)}\\
&\le &2^k {\delta(w_n)+c_0|w_n-z_n|\over \tau_{j,n}} {\delta(w_j) +c_0|w_j-z_j| \over \tau_{j, n} } \prod_{s=1}^\ell \prod_{i\in I_s}  {\delta(w_{i})+c_0|w_{i}-z_{i}|\over \tau_{i_s^*, i}(z,w)}\\
&\le &2^k {\delta(w_n)+c_0|w_n-z_n|\over \tau_{j,n}}  \prod_{s=1}^{\ell +1} \prod_{i\in I_s'}  {\delta(w_{i})+c_0|w_{i}-z_{i}|\over \tau_{i_s^*, i}(z,w)}\\
\end{eqnarray*}
with $I_1'=\{j\}$ and $I'_s=I_{s-1}$ for $2\le s\le k+1$. Therefore (\ref{S1}) holds for new $I'=I \cup\{ j \}$. Therefore, we have proved that (\ref{S}) and (\ref{S1}) hold for  $|I|=k+1$.  By the
Principle of Mathemtical Induction, we have proved (\ref{S}) and (\ref{S1}) hold.
Therefore, the proof of the theorem is complete.\epf

\begin{lemma} \label{L1} Let $0<\alpha<2$ and  $0<\beta<1$ and $2-\beta-\alpha > 0$. Then
$$
\int_\Omega {1\over |w_i-z_i|^\alpha} \delta(w_i)^{-\beta} dA(w_i)
\le C_\Omega C_{\alpha, \beta},
$$
where
$$
C_{\alpha,\beta} =\left\{ \begin{matrix} {1\over (1-\alpha)(1-\beta) },\quad \hbox{if } \ 0<\alpha, \beta<1;\cr {1\over (\alpha-1)(2-\alpha-\beta)},\quad \quad
\hbox{if } \ 1<\alpha<2-\beta;\cr
 {1\over (1-\beta)^2},\quad \hbox{ if } \ \alpha=1.\cr\end{matrix}\right.
$$
\end{lemma}

\proof Notice that for $0<\alpha<2$ and $2-\beta-\alpha> 0$, one has
\begin{eqnarray*}
\int_{D(z_i, \delta(z_i)/2)} |w_i-z_i|^{-\alpha} \delta(w_i)^{-\beta} dA(w_i)
&\le & 4^\beta \delta(z_i)^{-\beta}\int_{|w_i|<\delta(z_i)/2} |w_i|^{-\alpha} dA(w_i)\\
&\le & 2^{2\beta+1}\pi \delta(z_i)^{-\beta}\int_0 ^{\delta(z_i)/2} r^{-\alpha+1} dr\\
&=&{8\pi  \over 2-\alpha}\delta(z_i)^{2-\beta-\alpha}.
\end{eqnarray*}
Let $A$ be the diameter of $\Omega$ and $B$ the length of $\d \Omega$. Then, without loss of generality, one may assume that $\alpha\ne 1$,
since the case can be proved similarly.
\begin{eqnarray*}
\lefteqn{\int_\Omega (|w_i-z_i|+\delta(z_i))^{-\alpha} \delta(w_i)^{-\beta} d A(w) }\\
&\le&
C\int_0^A\int_0^B {1\over |x+iy|^\alpha } y^{-\beta} dx d y\\
 &\le&
C\int_0^A\int_0^B {1\over (x+ y)^\alpha } y^{-\beta} dx d y\\
&=&I_{\alpha, \beta}.
\end{eqnarray*}

If $0<\alpha <1$, then
$$
I_{\alpha, \beta}\le {B^{1-\alpha} \over 1-\alpha} {A^{1-\beta}\over 1-\beta}.
$$

If $\alpha >1$, then 
\begin{eqnarray*}
I_{\alpha,\beta}&\le&
{C\over (\alpha-1)} \int_0^A  y^{1-\alpha } y^{-\beta}  d y\\
&\le&
{C A^{2-\alpha-\beta} \over |\alpha-1|(2-\alpha-\beta)}.
\end{eqnarray*}

If $\alpha=1$, then 
\begin{eqnarray*}
I_{\alpha,\beta}&\le&
C \int_0^A ( -\ln y )  y^{-\beta} d y\\
&\le& {-\ln A \over 1-\beta}+{A^{1-\beta} \over (1-\beta)^2}.
\end{eqnarray*}
Therefore,  if  $0<\alpha<2$, $0<\beta<1$ and $2-\beta-\alpha >0$, then
\begin{eqnarray*}
\lefteqn{\int_\Omega |w_i-z_i|^{-\alpha} \delta(w_i)^{-\beta} dA(w_i)}\\
&=&\int_{D(z_i, \delta(z_i)/2)}  |w_i-z_i|^{-\alpha} \delta(w_i)^{-\beta} dA(w_i)
+\int_{\Omega\setminus D(z_i, \delta(z_i)/2)} |w_i-z_i|^{-\alpha} \delta(w_i)^{-\beta} dA(w_i)\\
&\le& {4\over 2-\alpha} +\int_{\Omega\setminus D(z_i, \delta(z_i)/2), \delta(w_i)\le \delta(z_i) } |w_i-z_i|^{-\alpha}\delta(w_i)^{-\beta} dA(w_i)\\
&&+\int_{\Omega\setminus D(z_i, \delta(z_i)/2), \delta(w_i)> \delta(z_i) } |w_i-z_i|^{-\alpha}\delta(w_i)^{-\beta} dA(w_i)\\
&\le&{C_\Omega C_{\alpha,\beta}}.
\end{eqnarray*}
The proof of the lemma is complete.\epf 

\medskip

\begin{lemma} \label{L2}
For $0<\beta <1$, $\alpha>0$  and $\beta+\alpha<2$, we define
 $$
  {\mathcal T}[g](\zeta)=\int_\Omega { 1\over \delta(\lambda)^{\beta} 
   |\lambda-\zeta|^{\alpha} }g(\lambda) dA(\lambda),\quad \zeta\in \Omega.
   $$
  Then
  
  (i) For any $1<p\le \infty$ such that  $p'\beta<1$ and $p'(\alpha+\beta)<2$, then  
  $$
  \Big|{\mathcal T}(g)(\zeta)\Big|\le \Big(C_\Omega  C_{p'\alpha, p'\beta}\Big)^{1/p'} \|g\|_{L^p(\Omega)},\quad \zeta\in \Omega.
  $$
  
  (ii) If $1< p\le 2$ such that $p\beta<1$ and $\alpha+p\beta<2$, then 
  ${\mathcal T} $  is bounded on $ L^q(\Omega)
  $ for all $p\le q\le p'$ and
  $$
  \|{\mathcal T}\|_{L^q(\Omega)\to L^q(\Omega)}\le C_\Omega C_{\alpha, p\beta};
  $$
  
  (iii) ${\mathcal T}$ is bounded on $L^p(\Omega)$ with $\|{\mathcal T}\|_{L^p(\Omega)\to L^p(\Omega)}\le C_{\alpha,\beta} C_\Omega$ for all $1<p\le \infty$.
 \end{lemma}
 
 \proof  (i) For any $g\in L^p(\Omega)$, since $p'\beta<1$ and $p'(\alpha+\beta)<2$, by Lemma \ref{L1}, one has
 $$
 |{\mathcal T}[g](\zeta)|  \le\Big( \int_\Omega \Big( {1\over  \delta(\lambda)^{\beta}   
   |\lambda-\zeta|^{\alpha} }\Big)^{p'} dA(\lambda)\Big)^{1/p'}\|g\|_{L^p}
   \le  \Big( C_\Omega C_{p'\alpha, p'\beta}\Big)^{1/p'} \|g\|_{L^p}.
$$
This proves Part (i). 

For Part (ii),  define
  $$
  r(\lambda)=\delta(\lambda)^{-{\beta\over p'}}.
  $$
  By Lemma \ref{L1}, one has
  \begin{eqnarray*}
  \Big( \int_\Omega  { \delta(\lambda)^{-\beta}   \over 
   |\lambda-\zeta|^{\alpha} } r(\zeta)^{p'} dA(\zeta)\Big)^{1/p'}
   &=&\Big( \int_\Omega  { \delta(\zeta)^{-{\beta}}   \over 
   |\lambda-\zeta|^{\alpha} }  dA(\zeta)\Big)^{1/p'}\delta(\lambda)^{-\beta/p'}\\
  & \le & \Big(  C_\Omega C_{\alpha,\beta} \Big )^{1/p'} \, r(\lambda).
\end{eqnarray*}

  For any $1<p< 2$ such that $p\beta<1$ and $\alpha+p\beta<2$,  by Lemma \ref{L1}, one has
  \begin{eqnarray*}
 \Big( \int_\Omega  { \delta(\lambda)^{-\beta }  \over 
   |\lambda-\zeta|^{\alpha}} r(\lambda)^p  dA(\lambda)\Big)^{1/p}
 &  = &  \Big(\int_\Omega{ \delta(\lambda)^{-p \beta}  \over 
   |\lambda-\zeta|^\alpha }  dA(\lambda)\Big)^{1/p}\\
   & \le & \Big(C_\Omega C_{\alpha, p\beta} \Big )^{1/p} \\
    & \le & \Big({C_\Omega C_{\alpha, p\beta} }\Big )^{1/p} r(\zeta).
  \end{eqnarray*}
  By Schur's lemma, this implies that 
  $$
  {\mathcal T}[g]=\int_\Omega { \delta(\lambda)^{-\beta}   \over 
   |\lambda-\zeta|^{\alpha} }g(\lambda) dA(\lambda)
   $$
   is bounded on the both $L^p$ and $L^{p'}$ with norm bounded by
   $$
   C=\Big({C_\Omega C_{\alpha, p\beta}}\Big)^{1/p'} \Big(C_\Omega C_{\alpha, p\beta })\Big)^{1/p}=C_\Omega C_{\alpha, p\beta}
   $$
By the interpolation theorem of integral operator, one has that ${\mathcal T}$ is bounded on $L^q$ with norm $C$ for all $p\le q\le p'$. This proves Part (ii). 

Part (iii) is the combination of Part (i) and Part (ii). For any $p>1$, we choose $1<p_0<\min\{p, 2\}$ such that
$$ 
p_0 \beta<1 \quad\hbox{and}\quad p_0(\alpha+\beta)<2.
$$
By Part (ii), we have ${\mathcal T}$ is bounded on $L^q(\Omega)$ for $p_0\le q\le p_0'$ and 
$$
\|{\mathcal T}\|_{L^q(\Omega)\to L^q(\Omega)} \le C_{\alpha, p_0\beta} C_\Omega.
$$
By Part (i), we have
$$
\|{\mathcal T}\|_{L^{p_0'}(\Omega)\to L^\infty(\Omega)}\le C_{p_0\alpha, p_0\beta}C_\Omega
$$
Therefore, let $p_0\to 1^+$, we have $C_{p_0\alpha, p_0\beta}\to C_{\alpha, \beta}$ and ${\mathcal T}$ is bounded on $L^p(\Omega)$ for all $1<p\le \infty$. Moreover,
$$
\|{\mathcal T}\|_{L^p(\Omega)\to L^p(\Omega)}\le C_{\alpha, \beta} C_\Omega,\quad 1<p\le \infty.
$$
 Therefore, the proof of the lemma is complete.\epf

For $n>1$, we propose to prove the following theorem.

\begin{theorem} Let $\Omega$ be a bounded domain in $\CC$ with $\delta_D(\cdot)\in \hbox{Lip}(\Omega)$.
Let $T$ be the linear integral operator defined by (\ref{T}) and (\ref{Tj}). Then there a positive constant  $C_\Omega =C(\|\delta_\Omega\|_{Lip})$ depending only on 
$\Omega$ such that
$$
\|T\|_{L^p(\Omega^n)\to L^p(\Omega^n)}\le  C_n C_\Omega^n ,\quad\hbox{ for all } \quad 1< p\le \infty.
$$
\end{theorem}

\proof  Let $I=\{i_1,\cdots, i_k\}=\cup_{s=1}^\ell I_s$ with $I_s\cap I_t=\emptyset$ if $s\ne t$;  $j\not\in I$  and $J=\{j, I\}$. 
For each $v\in I$, we choose  $p_v \in (2, \infty)$  with $1/p_v+1/p_v'=1$ .  Let
$$
\epsilon_s =\sum_{i\in I_s} {1\over p_i } <1.
$$
%For a fixed $z\in \Omega^n$, write
%$$
%\Omega=\Omega_{i_\ell}\cup \Omega_{i_\ell}^c
%$$
%where
%$$
%\Omega_{i_\ell}=\{w\in \Omega: |w_{i_\ell}-z_{i_\ell}|<\delta(w_{i_\ell})\}\quad
%\hbox{and}\quad \Omega^c_{I_\ell}=\Omega\setminus \times \Omega_{i_\ell}.
%$$
%Write
%$$
%\Omega^{k+1} =\sum_{|I|+|J|=k+1}\Omega_I\times \Omega_J^c, \quad \hbox{and }\ \Omega_{i_1}\times\cdots \Omega_{i_p}
%$$
For each $1\le s\le \ell$, we write
$$
I_s=I_s^1\cup I_s^2.
$$

For any $i\in I$, if $p_i\ge 2$ and $1/p_i+1/p_i'=1$, one has
$$
\tau_{i, j}(z,w) \ge {1\over p_i} |w_i-z_i|^2\delta(w_i)+{1\over p_i'} |w_j-z_j|^2 \delta(w_j)\ge (|w_i-w_i|^2 \delta(w_i))^{1/p_i} (|w_j-z_j|^2\delta(w_j))^{1/p_i'}.
$$
By Theorem \ref{Tj}, one has
 \begin{eqnarray*}
\lefteqn{ |S^j_J (z,w) | }\\
& \le& 
 {2^{k-1} \over |w_j-z_j| } \prod_{s=1}^{\ell} \prod_{i\in I_s} {\delta(w_{i})+c_0|w_{i}-z_{i}|\over \tau_{i_s^*, i} }\\
 &\le& {(2c_0)^k \over |w_j-z_j|} \prod_{s=1}^\ell \sum_{I_s=I_s^1\cup I_s^2} \prod_{u\in I_s^1}  {\delta(w_u)\over \tau_{i^*_s, u}}
 \prod_{v\in I_s^2}  {| w_v-z_v| \over \tau_{i_s^*, v} } \\
   & \le&{ (2c_0)^k\over |w_j-z_j|} \prod_{s=1}^\ell  \sum_{I_s=I_s^1\cup I_s^2} \prod_{u\in I_s^1} { \delta(w_{i_s^*})^{-1/p_{u} } \delta(w_u)^{1-1/p'_u}\over
   |w_{i_s^*}-z_{i_s^*}|^{2/p_u}
  |w_u-z_u|^{2/p_u'} }\prod_{v\in I_s^2}
   { \delta(w_{i_s^*})^{-1/p_v}  | w_v-z_v|^{1-2/p_v'}\over   |w_{i_s^*}-z_{i_s^*}|^{2/p_v}\delta(w_v)^{1/p_v'} } \\
   & =&{(2c_0)^k \over |w_j-z_j|}\prod_{s=1}^\ell\Big({ \delta(w_{i_s^*})^{-1}\over  |w_{i_s^*}-z_{i_s^*}|^2}\Big)^{\epsilon_s}
  \sum_{I_s=I_s^1\cup I_s^2}  \prod_{u\in I_s^1} {  \delta(w_u)^{1-1/p'_u}\over
  |w_u-z_u|^{2/p_u'} }\prod_{v\in I_s^2}
   {| w_v-z_v|^{1-2/p_v'} \over  \delta(w_v)^{1/p_v'} }.
 \end{eqnarray*}
 
 For each $1\le s\le \ell$, we let
 $$
 p_i=m_s>2,\quad \hbox{for any } i\in I_s,\quad \epsilon_s=\sum_{i\in I_s} {1\over p_i}={k_s\over m_s} \quad\hbox{with } k_s=|I_s|.
 $$
 Let $\epsilon>0$ be very small to be determined.  We choose $m_s>2$ such that
 $$
 \epsilon_s+{{1\over m'_{s-1}}}=1-k_s \cdots k_\ell \epsilon \iff \epsilon_s+1=1-k_s\cdots k_\ell \epsilon+{1\over m_{s-1}}\iff \epsilon_s+
 k_s\cdots k_\ell\epsilon={1\over m_{s-1}}.
  $$
 Take 
 $$
 \epsilon_\ell= k_\ell\epsilon>0
 $$
 Then 
 $$
{1\over  m_{\ell-1}}=2 k_\ell \epsilon,\quad \epsilon_{\ell-1}={k_{\ell-1}\over m_{\ell-1}}=k_{\ell-1} 2 k_\ell \epsilon \quad\hbox{and }\  {1\over m_{\ell-2}}
 =3k_{\ell-1} k_\ell\epsilon
 $$
 Thus,
 $$
 \epsilon_{\ell-2}={k_{\ell-2}\over m_{\ell-2}}=3 k_{\ell-2}k_{\ell-1}k_\ell \epsilon\quad\hbox{and}\quad {1\over m_{\ell-3}}
 =4 k_{\ell-2}k_{\ell-1} k_\ell \epsilon.
 $$
 Therefore,
 $$
 {1\over m_s}=(\ell+1-s) k_{s+1}\cdots k_\ell \epsilon,\quad \epsilon_1={k_1\over m_1}=\ell k_1\cdots k_\ell\epsilon.
 $$
 We choose ${1\over 2(n+1)^{n+1} }<\epsilon\le {1\over (n+1)^{n+1}}$ such that
 $$
 \epsilon_1\le 1-\epsilon.
 $$
 Therefore,
 \begin{equation} \label{eps1}
 \epsilon_s+{1\over m_{s-1}'}
 =1 - k_s\cdots k_\ell \epsilon \le 1-\epsilon,
 \end{equation}
 
\begin{equation}\label{eps2}
 2\epsilon_{s}+2/m'_{s-1}-1=1 -2k_s\cdots k_\ell \epsilon \in (1-{1\over n+1}, 1-{1\over (n+1)^{n+1}})
 \end{equation}
 and
 \begin{equation}\label{eps3}
 {3\epsilon_{s}+3/m'_{s-1}-1\le 2-3\epsilon.}
  \end{equation}
 
 Notice that $i_s^*\in I_{s-1}$. If   $ i_s^* =u\in I_{s-1}^1$, by (\ref{eps1}),we have
 $$
 \Big({ \delta(w_{i_s^*})^{-1}\over  |w_{i_s^*}-z_{i_s^*}|^2}\Big)^{\epsilon_s}
  {  \delta(w_u)^{1-1/m'_{s-1}}\over
  |w_u-z_u|^{2/m_{s-1}'} } = {  \delta(w_u)^{1-1/m'_{s-1} -\epsilon_s}\over
  |w_u-z_u|^{ 2\epsilon_s+ 2/m_{s-1}'} } \le {C_\Omega\over  |w_u-z_u|^{ 2-\epsilon} } 
  $$
  By  Lemma 4.1 in \cite{Li}, one has that the linear operator
  $$
  {\mathcal T}(g)(z_u)=\int_\Omega {  \delta(w_u)^{1-1/m'_{s-1} -\epsilon_s}\over
  |w_u-z_u|^{ 2\epsilon_s+ 2/m_{s-1}'} } g(w_u) dA(w_u)
$$
is bounded on $L^p(\Omega)$ for any $1\le p\le \infty$ and $\|{\mathcal T}\|_{L^p\to L^p}\le {C\over \epsilon} C_\Omega=C_n C_\Omega$.

If $i_s^*=v \in I_{s-1}^2$, then
 $$
 \Big({ \delta(w_{i_s^*})^{-1}\over  |w_{i_s^*}-z_{i_s^*}|^2}\Big)^{\epsilon_s}
   {| w_v-z_v|^{1-2/ m_{s-1}'} \over  \delta(w_v)^{1/m_{s-1}'} }= {| w_v-z_v|^{1-2/m_{s-1}' -2\epsilon_s} \over  \delta(w_v)^{\epsilon_s+ 1/m_{s-1}'} }
$$
By (\ref{eps1}), (\ref{eps2}), (\ref{eps3}) and Lemma \ref{L2}, one has that the linear operator
$$
{\mathcal T}(g)(z_v)=\int_{\Omega}  {| w_v-z_v|^{1-2/m_{s-1}' -2\epsilon_s} \over  \delta(w_v)^{\epsilon_s+ 1/m_{s-1}'} } g(w_v) dA(w_v)
$$
is bounded on $L^p(\Omega)$ for all $1<p\le \infty$ and $\|{\mathcal T}\|_{L^p\to L^p}\le C_{1-\epsilon, 1- 2n^n \epsilon} C_\Omega =C_n C_\Omega$.

This implies that
$$
S^j_{j, I}(f)(z) =\int_{\Omega^{k+1}} S^j_{j, I} (z,w) f_j(w) dv_{k+1}(w)
$$
is bounded on $L^p(\Omega^n )$ for all $1<p\le \infty$. Therefore, $S^j_{j,I}$ is bounded on $L^p(\Omega^n)$ for 
all $1<p\le \infty$ with $\|S^j_{j, I}\|_{L^p(\Omega^n)\to L^p(\Omega^n)}\le (C_n C_\Omega)^{n}$.

Combining all the above and the estimation formula (\ref{Tj})  in Theorem \ref{3.1}, one can easily see that $T_j$ is bounded on $L^p(\Omega^{n})$ for all $1<p\le \infty$ 
 with $\|T_j\|_{L^p(\Omega^n)\to L^p(\Omega^n) }\le ( C_n C_\Omega )^n$. 
 Apply Theorem \ref{3.1}, one has proved that $T$ is bounded on $L^p(\Omega^n)$ for all $1<p\le \infty$ with 
$$
\|T\|_{L^p(\Omega^n)\to L^p(\Omega^n)}\le  (C_nC_\Omega )^n.
$$
The proof of the theorem is complete. \epf

\section{The Proof of Theorem 1.1}

In this section, we will prove Theorem 1.1 by using the idea of the approximation which has been used in \cite{Li} and similar ideas
was also used in \cite{DPZ} and \cite{Yuan}.

Let $\Omega$ be a bounded domain in $\CC$ with the distance function $\delta_\Omega\in \hbox{Lip}(\Omega)$.
 For each $m\in \NN$, we define
\begin{equation}
\Omega_m=\{z\in \Omega: \delta(z)>1/m\}\textcolor{red}{.}
\end{equation}
Then $\Omega_m\subset \subset \Omega$  and $\lim_{m\to\infty}\Omega_{m}=\Omega$.  

For  any $\dbar$-closed  $f=\sum_{j=1}^{n}f_{j}d\zbar_{j}\in L^p_{(0,1)}(\Omega^n)$ and 
any $0<\epsilon<1/m$, we let $\chi\in C^\infty_0(D(0,1))$ be the cutoff function with $\int_\CC \chi(z) dA(z)=1$. We define
\begin{equation}
\chi_\epsilon(z)={1\over \epsilon^{2n}}\chi(z_1/\epsilon)\cdots \chi(z_n/\epsilon)
\end{equation}
and
\begin{equation}
f^\epsilon_{j}(z)=\chi_\epsilon * f_{j}(z)=\int_{\Omega^n} \chi_{\epsilon}(w) f_{j}(z-w) dv(w).
\end{equation}
Then  $f^\epsilon=\sum_{j=1}^{n}f^{\epsilon}_{j}d\zbar_{j}\in C^\infty_{(0,1)}(\overline{\Omega}_m^n )$ is $\dbar$-closed and $f^\epsilon\to f$ in $L^p(\Omega_m^n)$ as $\epsilon\to 0^+$ for any $m\in \NN$
and $1<p<\infty$. 

By the expression of $T_j$ in Theorem \ref{3.1},  for any $1<p<\infty$, one has 
$$
\|T_j[f_j^\epsilon ]-T_j[f_j]\|_{L^p(\Omega^n_m)}\le \Big( {C_n  C_ \Omega} \Big)^n \|f_j^\epsilon-f_j\|_{L^p(\Omega^n)} \to 0 \ \hbox{ as } \  \epsilon \to 0^+.
$$
This coupling with $\dbar \sum_{j=1}^n T_j[f_j^\epsilon ] =f^\epsilon $ on $\Omega_m$ implies that
$$  \dbar \sum_{j=1}^n T_j[f_j]=f\quad\hbox{on} \quad \Omega_m^n
$$
and
$$
\|\sum_{j=1}^n T_j[f_j]\|_{L^p(\Omega_m^n)}\le \Big( C_n  C_ \Omega \Big)^n ,\quad 1< p<\infty
$$
for any $m\in \NN$. Let $m\to +\infty$. Then
$$
\dbar \sum_{j=1}^n T_j[f_j]=f\quad \hbox{on} \quad \Omega^n
$$
and 
$$
\Big\|\sum_{j=1}^n T_j[f_j] \Big\|_{L^p(\Omega^n)} \le \Big( C_n C_ \Omega\Big)^n ,\quad 1< p<\infty.
$$
Let $p\to +\infty$. Then 
$$
\|T[f]\|_{L^\infty (\Omega^n)}=\Big\|\sum_{j=1}^n T_j[f_j] \Big\|_{L^\infty (\Omega^n)}\le  \Big( C_n C_\Omega\Big)^n.
$$
Therefore,
\begin{equation}
\|T[f]\|_{L^p(\Omega^n)}=\Big\|\sum_{j=1}^n T_j[f_j] \Big\|_{L^p(\Omega^n)}\le  \Big({ C_n}C_\Omega\Big)^n,\quad 1< p\le\infty
\end{equation}
and the proof of Theorem 1.1  is complete.\epf
\section*{Acknowledgement}
This article is partially supported by National Natural Science Foundation of China (Grant No. 12001259) and by Natural Science Foundation of Fujian Province (Grant No. 2020J01846).

\fontsize{11}{11}\selectfont
 \bigskip

\noindent Department of Mathematics, University of California, Irvine, CA 92697-3875, USA

\noindent Email addresses:  \ sli@uci.edu 
 \bigskip

 \noindent School of Mathematics and Data Science, Minjiang University, Fuzhou 350108, P. R. China

 \noindent Email addresses: \  lsj\_math@163.com
  \bigskip

 \noindent School of Mathematics and Statistics, FJKLMAA, Fujian Normal University, Fuzhou 350117,  P. R. China.

\noindent Email addresses:  \ jluo@fjnu.edu.cn
\end{document}